\documentclass[11pt]{article}
\usepackage{epsfig}
\usepackage{pst-grad} % For gradients
\usepackage{pst-plot} % For axes
\usepackage{amsmath,amssymb,amsthm,graphics,amscd, pst-plot,pstricks,mathrsfs}
\usepackage{pst-eps}
\usepackage{texdraw}

\newtheorem{theorem}{Theorem}
\newtheorem{lemma}{Lemma}

\def\square{\hbox{\vrule height8pt depth0pt
\vbox{\hrule width7.2pt\vskip7.2pt\hrule width7.2pt}\vrule
height8pt depth0pt}\smallskip}

\title{\bf Ramsey numbers of paths and graphs\\ of the same order\thanks{Supported by NSFC.}}
\author{Chaoping Pei, \hspace{3mm} Yusheng Li  \vspace{3mm} \\
{Department of Mathematics, Tongji University } \\
{Shanghai 200092, China} \\
{\small\em peichaoping@msn.cn, li\_yusheng@tongji.edu.cn}      }

\begin{document}

\date{}
\maketitle

\begin{abstract}

\indent

For graphs $F_n$ and $G_n$ of order $n$, if $R(F_n, G_n)=(\chi(G_n)-1)(n-1)+\sigma(G_n)$, then $F_n$ is said to be $G_n$-good, where $\sigma(G_n)$ is the minimum size of a color class among all proper vertex-colorings of $G_n$ with $\chi(G_n)$ colors.  Given $\Delta(G_n)\le \Delta$, it is shown that $P_n$ is asymptotically $G_n$-good if $\alpha(G_n)\le\frac{n}{4}$.

\medskip

\end{abstract}

\section{Introduction}

\indent

Let $F$ and $G$ be graphs. The Ramsey number $R(F,G)$ is defined to be the smallest $N$ such that any red-blue edge-coloring of $K_N$ contains either a red $F$ or a blue $G$. We shall write $|G|$, $\chi(G)$, $\Delta(G)$ and $\delta(G)$  as the order of $G$, the chromatic number, the maximum and minimum degrees of $G$, respectively. Denote by $\sigma(G)$  the minimum size of a color class among all proper vertex-colorings of $G$ with $\chi(G)$ colors. Burr \cite{burr} had the following general bound.

\begin{lemma}\label{burr} For any graph $G$, if $F$ is a connected graph with $|F|\ge \sigma(G)$, then
\begin{equation}\label{le-burr}
R(F,G)\ge (\chi(G)-1)(|F|-1)+\sigma(G).
\end{equation}
\end{lemma}
A connected graph $F$ is said to be $G$-good if the equality in $(\ref{le-burr})$ holds, in which a $K_p$-good graph is said to be $p$-good in short. Chv\'{a}tal \cite{chvatal} showed that a tree is $p$-good for any $p$. A family $\mathcal{F}$ of graphs is said to be $G$-good if all large graphs in $\mathcal{F}$ are $G$-good. Let $P_n^k$ be the $k$th power of $P_n$, whose edges consist of pairs $\{x,y\}$ with the distance in $P_n$ at most $k$. Let $bw(F)$ be the bandwidth of $F$ of order $n$, which is the smallest integer $k$ such that $F$ is a subgraph of $P_n^k$. Burr and Erd\H{o}s \cite{burr-edos} showed that for any $k$, the family of connected graphs $F$ with $bw(F)\le k$ is $p$-good for all $p$. Moreover, Nikiforov and Rousseau \cite{niki} showed that some larger families are $p$-good.

Burr \cite{burr} proved that, for any fixed connected graph $F$, the class of graphs homeomorphic to $F$ is {\em always-good}, i.e., $G$-good for any fixed $G$. Using Blow-up Lemma \cite{kom}, Allen, Brightwell and Skokan \cite{allen-bri} proved that the family of connected graphs $F$ with $\Delta(F)\le \Delta$ and $bw(F)=o(|F|)$ is always-good, where $\Delta$ is fixed.

Let us turn to $R(F_n,G_n)$, where $F_n$ and $G_n$ have the same order $n$. Burr \cite{burr} asked that if $\Delta(G_n)$ is bounded and $n$ is large, does
\[
R(G_n,G_n)=(\chi(G_n))-1)(n-1)+\sigma(G_n)?
\]
If the above equation holds, we say that $G_n$ is {\em itself-good}. Thus  $P_n$ is itself-good shown by Gerencs\'{e}r and Gy\'{a}rf\'{a}s \cite{ger}, and $C_n$ is itself-good shown by Faudree and Schelp \cite{fau}, and Rosta \cite{roa}. However, it is difficult to estimate $R(G_n,G_n)$ for general $G_n$ even with $\Delta(G_n)$ bounded, see \cite{allen-bri,chv-rod,fox}.

Note that the general answer for Burr's question is negative as shown in \cite{allen-bri} that $R(P_n^k,P_n^k)\ge (k+1)n-2k$. Recently, Pokrovskiy \cite{pok} proved $R(P_n,P_n^k)=(n-1)k+\lfloor \frac{n}{k+1}\rfloor$, which solves a conjecture in \cite{allen-bri}. For general $G_n$, it is difficult to answer whether $P_n$ is $G_n$-good or not, so we ask the question in a weaker form: if $\Delta(G_n)$ is bounded and $n\to \infty$, does
\[
R(P_n,G_n)=(\chi(G_n)-1)(n-1)+\sigma(G_n)+o(n)?
\]

The answer for the above question is positive in most cases. Let $\alpha(G_n)$ be the independence number of $G_n$, which is at least the largest size of a color class in any proper vertex-coloring of $G_n$. It is easy to see that $\sigma(G_n)\le \frac{n}{\chi(G_n)}$ and $\alpha(G_n)\ge \frac{n}{\chi(G_n)}\ge \frac{n}{\Delta(G_n)+1}$ for any graph $G_n$. A well-known fact \cite{alon-spe} states that  $\alpha(G_n)= \frac{n}{\Delta(G_n)+1}$ if and only if $G_n$ is a union of cliques of the same order, and thus we assume $\Delta\ge 3$ in the following result to avoid trivial cases.

\begin{theorem}\label{th-1}
Let $\Delta\ge 3$ be fixed, and $G_n$ a graph of order $n$ with $\Delta(G_n)\le \Delta$. If $n\to \infty$ and $\alpha(G_n)\le \frac{n}{4}$, then
\[
R(P_n,G_n)=(\chi(G_n)-1)n+\sigma(G_n)+o(n).
\]
\end{theorem}

\section{Lemmas}

In the following context, we always assume that $G_n$ is a graph of order $n$ with $\Delta(G_n)\le \Delta$.
To simplify notations, we write $k=\chi(G_n)$ and $\sigma=\sigma(G_n)$ instead of $k_n$ and $\sigma_n$, respectively. Let integers $\Delta$, $n$, $N$ and $\beta$, and real $\epsilon$ satisfy that
\[
\Delta\ge 3,  \hspace{3mm}  0<\epsilon\le \frac{1}{\Delta^5}, \hspace{3mm} \beta= \left\lceil \frac{1}{\epsilon^7}\right\rceil,\hspace{3mm}
N = (k-1)n+\sigma + \lceil\Delta^4\epsilon n\rceil,
\]
in which $n$ is sufficiently large for such given $\Delta$ and $\epsilon$.

For any red-blue edge-coloring of  $K_N$ on vertex set $V$, let $R$ and $B$ be the red and blue subgraph of the edge-colored $K_N$ respectively, where $\overline{R}=B$. We shall prove that if $R$ contains no red $P_n$, hence no $C_m$ with $m>n$, then $B$ contains $G_n$.

We will follow the method of \cite{allen-bri} to embed $G_n$ into the blue graph $B$. For some integer $r$, we call cycles
\[
C^{(1)},C^{(2)},\dots, C^{(r)}
\]
the {\em ordered disjoint longest cycles} in $R$, if $C^{(i)}$ is a longest cycle in $R\setminus \cup_{j=1}^{i-1}C^{(j)}$, $1\le i\le r$.
Denote by $c_i$ be the length of $C^{(i)}$. Now we assume that $r$ is the largest integer such that  $c_r\ge \epsilon^2n$. Then $r\le \frac{N}{\lceil \epsilon^2 n\rceil}<\frac{k}{\epsilon^2}$.  The following result is the Erd\H{o}s-Gallai theorem for cycles \cite{erdos-gallai}, in which $e(H)$ is the number of edges of $H$.
\begin{lemma}\label{erdos}
Let $H$ be a graph, and $c$ an integer with $3\le c\le |H|$. Then either $H$ contains a cycle of length at least $c$ or
\[
e(H)<\frac{(c-1)(|H|-1)}{2}+1.
\]
\end{lemma}
We shall define an often used set $W$ as follows.
\[
W =  V\setminus \cup_{j=1}^r V(C^{(j)}).
\]
\begin{lemma}\label{le-1-before}
There is a subset $W'$ of $W$ with $|W'|\ge (1-\epsilon) |W|$ such that every $\Delta$ vertices of $W'$ have at least $|W'|-\Delta\epsilon n$ common blue-neighbors in $W'$.
\end{lemma}
{\bf Proof.} Let $H$ be the subgraph of $R$ induced by $W$. Lemma \ref{erdos} with $c=\lceil\epsilon^2 n\rceil$ implies $e(H)< \epsilon^2n|H|/2$. Let $W_0\subseteq W$ such that any vertex in $W_0$ has at least $\epsilon n$ red-neighbors in $W$. Then $|W_0|\le \epsilon |W|$. Let $W'=W\setminus W_0$. Then $|W'|\ge (1-\epsilon) |W|$ and any vertex in $W'$ has at least $|W|-\epsilon n$ blue-neighbors in $W$. Hence every $\Delta$ vertices in $W'$ have at least $|W'|-\Delta\epsilon n$ common blue-neighbors in $W'$.  \hfill\square

\medskip

If $|W|\ge n+2\Delta\epsilon n$, then $|W'|>(1+\Delta\epsilon) n$, and any $\Delta$ vertices of $W'$  in Lemma \ref{le-1-before} have at least $n$ common blue-neighbors in $W'$, and thus the subgraph of $B$ induced by $W'$ contains $G_n$. Therefore in the following proof we assume that
\[
|W|< n+2\Delta\epsilon\, n, \hspace{3mm}\mbox{hence}\hspace{3mm} r\ge k-1
\]
as $\frac{N-|W|}{n}> \frac{(k-2)n+\sigma(G)}{n}>k-2$. This assumption will be used in the proof of Theorem \ref{th-1}.
Let $K_{n_1,\dots,n_r}$ be the complete $r$-partite graph with $n_i$ vertices in the $i$th part.  Note that
$n> c_1\ge \cdots \ge c_r\ge \epsilon^2 n$. The integer $\beta=\lceil \frac{1}{\epsilon^7}\rceil$ in the following result can be smaller, we assume so for further propose.

\begin{lemma}\label{le-2}
The blue graph $B$ contains a $K_{c_1-\beta,\dots,c_r-\beta}$ on partition $(V_1,\dots,V_r)$ with $V_i\subseteq V(C^{(i)})$ for $1\le i\le r$.
\end{lemma}
{\bf Proof.} For any $C^{(i)}$ and $C^{(j)}$ with $1\le i<j\le r$, we {\bf claim} that there exist $S\subseteq C^{(i)}$ and $T\subseteq C^{(j)}$ with $|S|\ge c_i-\frac{9}{\epsilon^4}$  and $|T|\ge c_j- \frac{9}{\epsilon^4}$ such that $S$ and $T$ are completely adjacent in $B$.
In fact, we partition each $C^{(\ell)}$ into segments of consecutive vertices, in which each segment contains $\lceil\frac{c_r}{2}\rceil$ vertices but at most one contains less. Let $S'$ be a segment of $C^{(i)}$, and $T'$ of $C^{(j)}$. Clearly, there are no two independent edges between $S'$ and $T'$ in $R$, since otherwise, we can enlarge $C^{(i)}$ by adding the independent edges and edges of $C^{(j)}\setminus T'$ as illustrated in Figure 1, impossible. Thus we can ignore at most one vertex in each of $S'$ and $T'$ such that the remaining vertices in $S'$ and that in $T'$ are completely adjacent in $B$.

\medskip

\begin{center}
\scalebox{0.8} % Change this value to rescale the drawing.
{
\begin{pspicture}(0,-1.9084375)(11.36,1.8084375)
\definecolor{color6}{rgb}{0.6,0.6,0.6}
\usefont{T1}{ptm}{m}{n}
\rput(2.0403125,0.5584375){\small $S'$}
\usefont{T1}{ptm}{m}{n}
\rput(4.2803125,0.6184375){\small $T'$}
\usefont{T1}{ptm}{m}{n}
\rput(2.776406,-1.6015625){possible}
\usefont{T1}{ptm}{m}{n}
\rput(9.198594,-1.6815625){impossible}
\pscircle[linewidth=0.04,dimen=outer](1.32,0.6284375){1.12}
\pscircle[linewidth=0.04,dimen=outer](4.54,0.6484375){0.76}
\psarc[linewidth=0.12,linecolor=color6](1.2,0.6084375){1.2}{316.16913}{45.0}
\psarc[linewidth=0.12,linecolor=color6](4.55,0.6384375){0.75}{102.72436}{254.24883}
\psline[linewidth=0.04cm,dotsize=0.07055555cm 2.0]{*-*}(2.28,1.1284375)(3.78,0.6484375)
\psline[linewidth=0.04cm,dotsize=0.07055555cm 2.0]{*-*}(2.38,0.6284375)(3.78,0.6484375)
\psline[linewidth=0.04cm,dotsize=0.07055555cm 2.0]{*-*}(2.36,0.1884375)(3.76,0.6484375)
\usefont{T1}{ptm}{m}{n}
\rput(7.9603133,0.6984375){\small ${S'}$}
\usefont{T1}{ptm}{m}{n}
\rput(10.340312,0.6784375){\small ${T'}$}
\usefont{T1}{ptm}{m}{n}
\rput(3.6303124,0.9784375){\small ${v}$}
\pscircle[linewidth=0.04,dimen=outer](7.38,0.6884375){1.12}
\pscircle[linewidth=0.04,dimen=outer](10.6,0.7084375){0.76}
\psarc[linewidth=0.12,linecolor=color6](7.26,0.6684375){1.2}{316.16913}{45.0}
\psarc[linewidth=0.12,linecolor=color6](10.61,0.6984375){0.75}{102.72436}{254.24883}
\psline[linewidth=0.04cm,dotsize=0.07055555cm 2.0]{*-*}(8.3,1.2884375)(9.96,1.1084375)
\psline[linewidth=0.04cm,dotsize=0.07055555cm 2.0]{*-*}(8.3,0.0284375)(9.94,0.3484375)
\usefont{T1}{ptm}{m}{n}
\rput(1.2842188,-0.8015625){$C^{(i)}$}
\usefont{T1}{ptm}{m}{n}
\rput(7.3442187,-0.7815625){$C^{(i)}$}
\usefont{T1}{ptm}{m}{n}
\rput(4.3942184,-0.7615625){$C^{(j)}$}
\usefont{T1}{ptm}{m}{n}
\rput(10.514218,-0.7815625){$C^{(j)}$}
\end{pspicture}
}

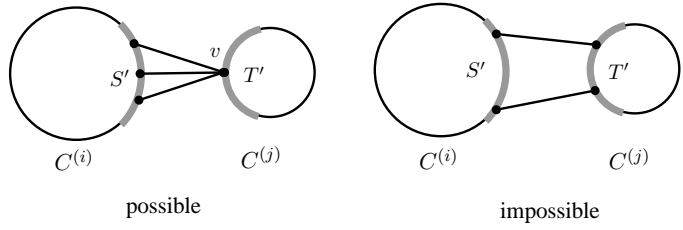
{\bf Figure 1} Possible and impossible edges between $C^{(i)}$ and $C^{(j)}$ in $R$
\end{center}

\medskip

Specifically for $C^{(i)}$ and $C^{(r)}$, the later contains two segments $T_1'$ and $T_2'$ of consecutive vertices with $|T'_2|\le |T_1'|=\lceil\frac{c_r}{2}\rceil$. By the mentioned process, we ignore one vertex in $T_1'$ when considering it with a segment of $C^{(i)}$. As $C^{(i)}$ has at most $\lceil\frac{c_j}{|T'|}\rceil\le\lceil\frac{2n}{c_r}\rceil \le \frac{3}{\epsilon^2}$ such segments, we ignore at most $\frac{3}{\epsilon^2}$ vertices in each of $T_1'$ and $C^{(i)}$ such that the remaining vertices in $T'_1$ and that in $C^{(i)}$ are completely adjacent in $B$. Same argument holds for $T_2'$ and $C^{(i)}$. By considering each of $C^{(1)},\dots,C^{(r-1)}$ with $C^{(r)}$, we can ignore at most $\frac{6}{\epsilon^2}$ vertices in each of $C^{(1)},\dots,C^{(r-1)}$ and at most $\frac{6r}{\epsilon^2}$ vertices in $C^{(r)}$ such that the remaining vertices in each of $C^{(1)},\dots,C^{(r-1)}$ and that in $C^{(r)}$ are completely adjacent in $B$.

Repeat the process to $C^{(1)},\dots,C^{(i-1)}$ and $C^{(i)}$ for each $1\le i\le r-1$, and thus we can ignore at most $\frac{9r}{\epsilon^4}\le \beta$ vertices in each $C^{(i)}$ such that the remaining vertices induce a subgraph of $B$ that contains $K_{c_1-\beta,\dots,c_r-\beta}$ as required.   \hfill\square

\medskip

The following  result says that the common blue-neighborhood of any $\Delta$ vertices of $W$ covers almost half vertices of each $C^{(i)}$, in which the idea of the proof comes from \cite{burr}.

\begin{lemma} \label{le-3}
Let $w_1,\dots,w_\Delta$ be vertices of $W$. Then the common blue-neighborhood $\cap_{j=1}^{\Delta}N_{B}(w_j)$ contains at least $\frac{1}{2}(c_i-\Delta^2)$ vertices of $C^{(i)}$ for $1\le i\le r$.
\end{lemma}
{\bf Proof.} Fix $i$ with $1\le i\le r$. Let vertices in $V(C^{(i)})$ be labeled by $\{1,2,\dots, c_i\}$ clockwise. For $1\le j\le \Delta$, let
\[
Z_j=N_R(w_j)\cap V(C^{(i)})= \{y_1,y_2,\dots \}, \hspace{3mm}\mbox{and}\hspace{3mm} Z_j+1=\{y_1+1,y_2+1,\dots \},
\]
where $Z_j+1$ contains vertices of $C^{(i)}$ next to that of $Z_j$. Clearly each vertex in $Z_j+1$ is non-adjacent to $w_j$ in $R$ from the maximality of $C^{(i)}$. For $\ell\neq j$ with $1\le \ell\le \Delta$, the vertex $w_\ell$ has at most one red-neighbor in $Z_j+1$, otherwise $C^{(i)}$ could be enlarged as illustrated in Figure 2.

\medskip

\begin{center}
\scalebox{0.7} % Change this value to rescale the drawing.
{
\begin{pspicture}(0,-2.33)(6.18,2.33)
\psellipse[linewidth=0.04,dimen=outer](3.09,-1.18)(3.09,1.15)
\psdots[dotsize=0.12](1.86,1.51)
\psdots[dotsize=0.12](4.64,1.53)
\psline[linewidth=0.04cm](1.82,1.51)(1.26,-0.25)
\psline[linewidth=0.04cm](1.9,1.47)(4.24,-0.15)
\psline[linewidth=0.04cm](4.62,1.49)(1.72,-0.15)
\psline[linewidth=0.04cm](4.62,1.49)(4.84,-0.23)
\psdots[dotsize=0.12](1.26,-0.23)
\rput(1.40,-0.6){$j_a$}\rput(2.20,-0.6){$j_a+1$}\rput(4.40,-0.6){$j_b$}\rput(5.20,-0.6){$j_b+1$}
\psdots[dotsize=0.12](1.76,-0.17)
\psdots[dotsize=0.12](4.18,-0.15)
\psdots[dotsize=0.12](4.84,-0.27)
\rput(1.79,1.92){$w_j$}
\rput(4.7,2.01){$w_\ell$}
\psframe[linewidth=0.04,linecolor=white,dimen=outer,fillstyle=solid](1.66,-0.11)(1.3,-0.27)
\psframe[linewidth=0.04,linecolor=white,dimen=outer,fillstyle=solid](4.76,0.01)(4.26,-0.33)
\rput(5.91,-0.26){$C^{(i)}$}
\end{pspicture}
}

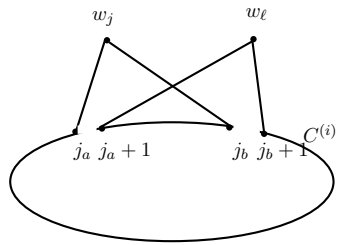
{\bf Figure 2} \hspace{3mm} a longer cycle
\end{center}

\medskip

After considering $w_1,\dots,w_{\Delta}$ in the same way, we know that $(Z_j+1)\cap(\cap_{\ell=1}^{\Delta}N_{B}(w_\ell))$ contains at least $|Z_j+1|-\Delta$ vertices. Inductively, we have
\[
\big|(\cup_{j=1}^\Delta(Z_j+1))\cap(\cap_{\ell=1}^{\Delta}N_{B}(w_\ell))\big |\ge |\cup_{j=1}^\Delta (Z_j+1)|-\Delta^2.
\]
Let $x=|\cap_{\ell=1}^{\Delta}N_{B}(w_\ell)\cap V(C^{(i)})|$, the number of common blue-neighbors of $w_1,\dots, w_{\Delta}$ in $V(C^{(i)})$. Then, we have $c_i-x=|\cup_{j=1}^\Delta Z_j|$ and
\[
x\ge |\cup_{j=1}^\Delta(Z_j+1)|-\Delta^2 = |\cup_{j=1}^\Delta Z_j|-\Delta^2,
\]
implying $x\ge \frac{c_i-\Delta^2}{2}$.\hfill\square

\medskip

Lemma \ref{le-2} says that the cycles $C^{(1)},C^{(2)}, \dots, C^{(r)}$ are almost completely connected each other in $B$, and Lemma \ref{le-3} tells us that edges between $W$ and $\cup_{\ell=1}^r C^{(\ell)}$ are dense in $B$.

\section{Proof of the main results}

Recall $W = V\setminus \cup_{j=1}^{r} V(C^{(j)})$ and $r\ge k-1$, and define
\[
U  = V\setminus \cup_{j=1}^{k-1} V(C^{(j)}).
\]
For disjoint subsets $Q$ and $Q'$ of $U$, we say that  $Q$ and $Q'$ {\em share none of} $C^{(k)},C^{(k+1)},\dots,C^{(r)}$, if at least one of $Q\cap V(C^{(i)})$ and $Q'\cap V(C^{(i)})$ is empty for any $k\le i\le r$. When $r=k-1$, any disjoint $Q$ and $Q'$ hold the property trivially. Let $\{1,2,\dots,k\}$ be the set of colors used for a proper vertex-coloring of $G_n$, and $A_i$ the color class of vertices in color $i$  such that $a_1\ge a_2\ge\dots\ge a_k=\sigma$, where $a_i=|A_i|$ for $1\le i\le k$. We shall embed $A_1,A_2,\dots,A_k$ into $B$, in which if we choose a vertex from some $C^{(i)}$ for a vertex of $A_a$, then we will not choose any from the same $C^{(i)}$ for $A_b$, where $a\neq b$.

By Lemma \ref{le-2}, there is a $K_{c_1-\beta,\dots,c_{r}-\beta}$ in $B$ with partition $(V_1,\dots, V_r)$ such that $V_i\subseteq V(C^{(i)})$ for $1\le i\le r$, and from Lemma \ref{le-1-before}, there is a subset $W'$ of $W$ with $|W'|\ge (1-\epsilon)|W|\ge|W|-2\epsilon\, n$ such that every $\Delta$ vertices in $W'$ have at least $|W'|-\Delta\epsilon\, n$ common blue-neighbors in $W'$. Let $W_0=W\setminus W'$. Then $|W_0|\le 2\epsilon\, n$. Let $\overline{c}$ be the length of the largest cycle in $U$. Clearly, $\overline{c}=c_k$ if $r\ge k$ and $\overline{c}<\epsilon^2 n$ otherwise.

\medskip

{\bf The Outline of the remaining proofs. } Partition  $U\setminus W$ or $U\setminus W_0$ into $Q_1,Q_2,\dots,Q_k$, where some $Q_i$ is possibly empty, such that any pair $Q_i$ and $Q_j$ share none of $C^{(k)}, \dots, C^{(r)}$, $A_k$ can be embedded into $Q_k$, and $A_i$ can be embedded into $V_i\cup Q_i$ for $1\le i\le k-1$. {\em We will embed $A_i$ into $V_i\cup Q_i$ by an embedding $\phi$ such that for every $v\in A_i$ and its neighbors $u_1,u_2,\dots$, $\phi(v)$ is adjacent to $\phi(u_1),\phi(u_2),\dots$ in $B$.} The outline is illustrated by Figure 3.

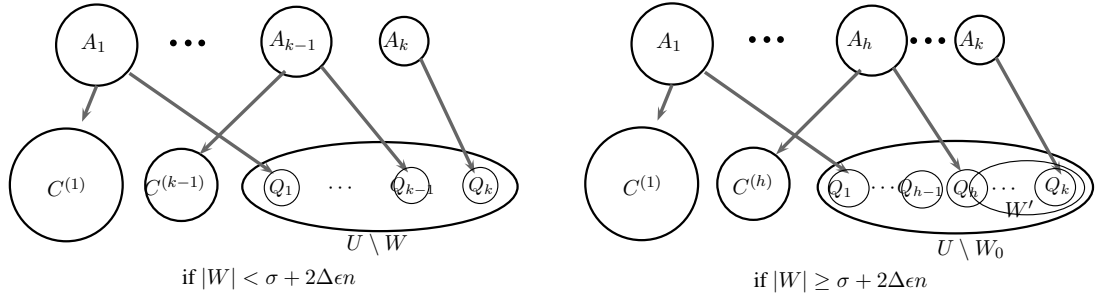
\begin{center}
\scalebox{0.75} % Change this value to rescale the drawing.
{
\begin{pspicture}(0,-2.6089063)(19.257032,2.5689063)
\definecolor{color130}{rgb}{0.4,0.4,0.4}
\pscircle[linewidth=0.04,dimen=outer](1.56,1.8089062){0.74}
\pscircle[linewidth=0.04,dimen=outer](5.1,1.8689063){0.64}
\psdots[dotsize=0.12](2.92,1.8489063)
\psdots[dotsize=0.12](3.18,1.8489063)
\psdots[dotsize=0.12](3.42,1.8489063)
\pscircle[linewidth=0.04,dimen=outer](7.01,1.8789061){0.45}
\pscircle[linewidth=0.04,dimen=outer](1.02,-0.67109376){1.02}
\psellipse[linewidth=0.04,dimen=outer](6.58,-0.7010938)(2.46,0.81)
\pscircle[linewidth=0.04,dimen=outer](3.06,-0.6710938){0.66}
\psline[linewidth=0.06cm,linecolor=color130,arrowsize=0.05291667cm 2.0,arrowlength=1.4,arrowinset=0.4]{<-}(1.3,0.44890618)(1.56,1.0889063)
\psline[linewidth=0.06cm,linecolor=color130,arrowsize=0.05291667cm 2.0,arrowlength=1.4,arrowinset=0.4]{<-}(4.7,-0.47109374)(2.14,1.3089062)
\psline[linewidth=0.06cm,linecolor=color130,arrowsize=0.05291667cm 2.0,arrowlength=1.4,arrowinset=0.4]{<-}(3.38,-0.13109376)(4.9,1.3489063)
\psline[linewidth=0.06cm,linecolor=color130,arrowsize=0.05291667cm 2.0,arrowlength=1.4,arrowinset=0.4]{<-}(7.08,-0.45109376)(5.54,1.4289062)
\psline[linewidth=0.06cm,linecolor=color130,arrowsize=0.05291667cm 2.0,arrowlength=1.4,arrowinset=0.4]{<-}(8.26,-0.47109374)(7.3,1.5689063)
\usefont{T1}{ptm}{m}{n}
\rput(1.505625,1.8389063){$A_1$}
\usefont{T1}{ptm}{m}{n}
\rput(5.0356255,1.8789061){$A_{k-1}$}
\usefont{T1}{ptm}{m}{n}
\rput(6.865625,1.8789061){$A_k$}
\usefont{T1}{ptm}{m}{n}
\rput(1.0356249,-0.7010938){$C^{(1)}$}
\usefont{T1}{ptm}{m}{n}
\rput(2.9456248,-0.6810938){$C^{(k-1)}$}
\usefont{T1}{ptm}{m}{n}
\rput(4.8317184,-0.72109383){\small ${Q_1}$}
\usefont{T1}{ptm}{m}{n}
\rput(7.161719,-0.7010938){\small ${Q_{k-1}}$}
\usefont{T1}{ptm}{m}{n}
\rput(8.371719,-0.7010938){\small ${Q_k}$}
\usefont{T1}{ptm}{m}{n}
\rput(6.5356255,-1.7410939){$U\setminus W$}
\pscircle[linewidth=0.04,dimen=outer](11.76,1.8289063){0.74}
\pscircle[linewidth=0.04,dimen=outer](15.3,1.8889061){0.64}
\pscircle[linewidth=0.04,dimen=outer](17.21,1.8989061){0.45}
\pscircle[linewidth=0.04,dimen=outer](11.22,-0.6510937){1.02}
\psellipse[linewidth=0.04,dimen=outer](16.78,-0.7110937)(2.46,0.84)
\psdots[dotsize=0.12](13.18,1.9089061)
\psdots[dotsize=0.12](13.44,1.9089061)
\psdots[dotsize=0.12](13.68,1.9089061)
\pscircle[linewidth=0.04,dimen=outer](13.2,-0.6510937){0.66}
\psline[linewidth=0.06cm,linecolor=color130,arrowsize=0.05291667cm 2.0,arrowlength=1.4,arrowinset=0.4]{<-}(11.5,0.4689062)(11.76,1.108906)
\psline[linewidth=0.06cm,linecolor=color130,arrowsize=0.05291667cm 2.0,arrowlength=1.4,arrowinset=0.4]{<-}(14.9,-0.47109374)(12.34,1.3289063)
\psline[linewidth=0.06cm,linecolor=color130,arrowsize=0.05291667cm 2.0,arrowlength=1.4,arrowinset=0.4]{<-}(13.6,-0.13109376)(15.1,1.3689063)
\psline[linewidth=0.06cm,linecolor=color130,arrowsize=0.05291667cm 2.0,arrowlength=1.4,arrowinset=0.4]{<-}(16.88,-0.43109375)(15.68,1.4089062)
\psline[linewidth=0.06cm,linecolor=color130,arrowsize=0.05291667cm 2.0,arrowlength=1.4,arrowinset=0.4]{<-}(18.66,-0.49109384)(17.5,1.5889063)
\usefont{T1}{ptm}{m}{n}
\rput(11.705625,1.858906){$A_1$}
\usefont{T1}{ptm}{m}{n}
\rput(15.085625,1.8989061){$A_{h}$}
\usefont{T1}{ptm}{m}{n}
\rput(17.065626,1.8989061){$A_k$}
\usefont{T1}{ptm}{m}{n}
\rput(11.235625,-0.6810939){$C^{(1)}$}
\usefont{T1}{ptm}{m}{n}
\rput(13.195625,-0.66109383){$C^{(h)}$}
\usefont{T1}{ptm}{m}{n}
\rput(14.731718,-0.7410938){\small ${Q_1}$}
\usefont{T1}{ptm}{m}{n}
\rput(16.97172,-0.7410939){\small ${Q_h}$}
\usefont{T1}{ptm}{m}{n}
\rput(18.59172,-0.7010938){\small ${Q_k}$}
\usefont{T1}{ptm}{m}{n}
\rput(17.055626,-1.8210938){$U\setminus W_0$}
\usefont{T1}{ptm}{m}{n}
\rput(4.616875,-2.341094){if $|W|< \sigma+2\Delta\epsilon n$}
\usefont{T1}{ptm}{m}{n}
\rput(14.746875,-2.3810937){if $|W|\ge \sigma+2\Delta\epsilon n$}
\psellipse[linewidth=0.02,dimen=outer](18.05,-0.7210938)(1.03,0.49)
\psellipse[linewidth=0.02,dimen=outer](4.84,-0.72109383)(0.32,0.33)
\psellipse[linewidth=0.02,dimen=outer](7.14,-0.6810939)(0.32,0.33)
\psellipse[linewidth=0.02,dimen=outer](8.36,-0.6810939)(0.32,0.33)
\usefont{T1}{ptm}{m}{n}
\rput(17.915627,-1.1010936){${W'}$}
\psellipse[linewidth=0.02,dimen=outer](18.55,-0.7110937)(0.37,0.34)
\psellipse[linewidth=0.02,dimen=outer](16.99,-0.7310938)(0.37,0.34)
\psellipse[linewidth=0.02,dimen=outer](14.89,-0.7310938)(0.37,0.34)
\usefont{T1}{ptm}{m}{n}
\rput(5.912812,-0.7210938){$\cdots$}
\usefont{T1}{ptm}{m}{n}
\rput(15.532812,-0.7210938){$\cdots$}
\psellipse[linewidth=0.02,dimen=outer](16.19,-0.7510938)(0.37,0.34)
\usefont{T1}{ptm}{m}{n}
\rput(16.161716,-0.7410939){\small ${Q_{h-1}}$}
\psdots[dotsize=0.12](16.04,1.8889061)
\psdots[dotsize=0.12](16.3,1.8889061)
\psdots[dotsize=0.12](16.54,1.8889061)
\usefont{T1}{ptm}{m}{n}
\rput(17.672812,-0.74109375){$\cdots$}
\end{pspicture}
}

{\bf Figure 3} \hspace{3mm} outline of the embedding
\end{center}

If $|W|< \sigma+2\Delta\epsilon n$, we partition $U\setminus W$ into $Q_1,\dots,Q_k$ such that any pair of $C^{(1)},\dots,C^{(k-1)},$ $Q_1,\dots,Q_k$ are almost completely connected in $B$. If $|W|\ge \sigma+2\Delta\epsilon n$, we partition $U\setminus W_0$ into $Q_1,\dots,Q_{h-1},Q_h,Q_{h+1},\dots,Q_k$ for some $h\le k-1$ such that $Q_1,\dots, Q_{h-1}\subseteq U\setminus W$, and $Q_{h+1},\dots,Q_k\subseteq W'$, and $Q_h$ consists of remaining vertices in $U\setminus W_0$. We shall embed $A_k$ into $Q_k$ first. Then by Lemma \ref{le-1-before}, for each $i$ with $h\le i\le k-1$, we shall embed a part of $A_i$ into $Q_i$ such that there are at most $\Delta^2\epsilon n$ vertices left in $W'$. Then, as Lemma \ref{le-1-before} and Lemma \ref{le-3} says that any $\Delta$ vertices in $W'$ have at least $(|V_i\cup Q_i|-r\Delta^2-\Delta\epsilon n)/2$ common blue-neighbors in $V_i\cup Q_i$, we shall put the remaining vertices in $A_i$ into $V_i\cup Q_i$, which will be explained in details.

Define
\[
\eta_i  =  \left\lfloor\frac{c_i-2\beta}{2}\right\rfloor
\]
for $1\le i\le r$.  Lemma \ref{le-2} and Lemma \ref{le-3} imply that any $\Delta$ vertices in $W$ have at least $\eta_i$ common blue-neighbors in $V_i$ as $\eta_i<\frac{1}{2}(c_i-\Delta^2-\beta)$. For $1\le i\le k-1$, let us denote by
\begin{eqnarray*}
\Lambda & = & \big\{j:\,1\le j\le k-1, a_j>\eta_j\big\},\\
\Gamma  & = & \big\{j:\,1\le j\le k-1,\, a_j>c_j-\beta\big\},
\end{eqnarray*}
and $\lambda=|\Lambda|$, $\gamma=|\Gamma|$. Then $\Gamma\subseteq \Lambda$. If $\Lambda=\emptyset$, the proof is trivial as $U$ is large enough for $A_k$ and $V_i$ is large enough for $A_i$ for $1\le i\le k-1$. So we assume that $\Lambda\neq \emptyset$.

\medskip

 Assume that $p\in \Lambda$. Then we have
\[\overline{c}\le c_p\le 2\eta_p+3\beta< 2a_p+3\beta\le 2\alpha(G_n)+3\beta\le n/2+3\beta.
\]
{\bf Case 1.} $ |W|<\sigma+2\Delta\epsilon n$. As $n\ge 2\alpha(G_n)+\overline{c}-3\beta\ge a_i+\sigma+\overline{c}-3\beta$ for $1\le i\le k-1$, we have
\begin{eqnarray*}
|U\setminus W|+\sum_{i\in\Gamma}|V_i| & \ge & N-(k-1-\gamma)n-(\sigma+2\Delta\epsilon n)-r\beta \\
& \ge & \gamma n+2\epsilon n\ge \sum_{i\in \Gamma}(a_i+\overline{c})+\sigma+\epsilon n.
\end{eqnarray*}
Partition $U\setminus W$ into $Q_1,\dots, Q_k$ such that any pair of them share none cycle of $C^{(k)},\dots,C^{(r)}$ as follows.
\begin{itemize}
\item For $i\not\in\Gamma\cup\{k\}$, let $Q_i=\emptyset$.
\item Note that $|V_i|< c_i< a_i+r\beta$ for $i\in\Gamma$ and $c_r\le \cdots\le c_k\le \overline{c}$, we can put a whole cycle among $C^{(k)},\dots,C^{(r)}$ one by one into $Q_i$ until $a_i+r\beta\le |V_i\cup Q_i|\le a_i+\overline{c}+r\beta$. Then for $Q_j$ with $j\in\Gamma\setminus\{i\}$ that has not been constructed yet, we put a whole unused cycle among $C^{(k)},\dots,C^{(r)}$ one by one into $Q_j$ similarly.
\item  Set $Q_k=(U\setminus W)\setminus \cup_{i\in\Gamma} Q_i$. Then $|Q_k|\ge \sigma+r\beta$.
\end{itemize}

Then the embedding can be constructed by Figure 3 and the explantation thereafter.

{\bf Case 2.} $|W|\ge \sigma+2\Delta\epsilon n$. As $n\ge 2\alpha(G_n)+\overline{c}-3\beta\ge 2a_i+\overline{c}-3\beta$ for $1\le i\le k-1$, we have
\begin{eqnarray*}
|U\setminus W_0|+\sum_{i\in\Lambda}|V_i|\ge N-(k-1-\lambda)n\ge \sum_{i\in \Lambda}(2a_i+\overline{c}+2\Delta\epsilon n)+\sigma.
\end{eqnarray*}

Partition $U\setminus W_0$ into $Q_1,\dots,Q_h,\dots,Q_k$ with $1\le h\le k-1$ such that $Q_1,\dots, Q_{h-1}\subseteq U\setminus W$, $Q_{h+1},\dots,Q_k\subseteq W'$, and any pair $Q_i$ and $Q_j$ share none cycle of $C^{(k)},\dots,C^{(r)}$ as follows.
\begin{itemize}
\item For $i\not\in \Lambda\cup\{k\}$, let $Q_i=\emptyset$.
\item Choose any $\sigma+\Delta\epsilon n$ vertices in $W'$ as $Q_k$.
\item  Let $h$ be the smallest integer such that
\[
|W'|\ge \sum_{j\in \Lambda \atop
h<j<k} (2a_j-|V_j|)+\sigma+(k-h)\Delta\epsilon n,
\]
and if $h=k-1$, then $|W'|\ge \sigma+\Delta\epsilon n$. For each $i\in\Lambda\cap\{h+1,\dots,k-1\}$, choose any $2a_i-|V_i|+\Delta\epsilon n$ vertices in $W'$ as $Q_i$.
\item For each $i\in \Lambda\cap \{1,\dots,h-1\}$, we put a whole cycle among $C^{(k)},\dots,C^{(r)}$ one by one into $Q_i$ until $2a_i+2\Delta\epsilon n\le |V_i\cup Q_i|\le 2a_i+\overline{c}+2\Delta\epsilon n$. We can do so since $c_i<2\eta_i+3\beta< 2a_i+2\Delta\epsilon n$ for $i\in\Lambda$ and $c_r\le \cdots\le c_k\le \overline{c}$.
\item  Let $Q_h=(U\setminus W_0)\setminus \cup_{i\neq h} Q_i$. Then $|V_h\cup Q_h|\ge 2a_h+2\Delta\epsilon n$.
\end{itemize}
Then the embedding can be constructed by Figure 3 and the explantation thereafter.\hfill\square

\indent

\end{document}